\documentclass[a4paper,11pt,leqno,twoside]{article}

\usepackage{amsmath, amsthm, amssymb}
\usepackage{graphics}

\makeatletter
\newif\ify@autoscale \y@autoscaletrue \def\Yautoscale#1{\ifnum #1=0
  \y@autoscalefalse\else\y@autoscaletrue\fi}
\newdimen\y@b@xdim
\newdimen\y@boxdim \y@boxdim=13pt
\def\Yboxdim#1{\y@autoscalefalse\y@boxdim=#1}
\newdimen\y@linethick    \y@linethick=.3pt
\def\Ylinethick#1{\y@linethick=#1}
\newskip\y@interspace \y@interspace=0ex plus 0.3ex
\def\Yinterspace#1{\y@interspace=#1}
\newif\ify@vcenter   \y@vcenterfalse
\def\Yvcentermath#1{\ifnum #1=0 \y@vcenterfalse\else\y@vcentertrue\fi}
\newif\ify@stdtext   \y@stdtextfalse
\def\Ystdtext#1{\ifnum #1=0 \y@stdtextfalse\else\y@stdtexttrue\fi}
\newif\ify@enable@skew   \y@enable@skewfalse
\y@vcentertrue
\def\y@vr{\vrule height0.8\y@b@xdim width\y@linethick depth 0.2\y@b@xdim}
\def\y@emptybox{\y@vr\hbox to \y@b@xdim{\hfil}}
\ify@enable@skew
 \def\y@abcbox#1{\if :#1\else
   \y@vr\hbox to \y@b@xdim{\hfil#1\hfil}\fi}
 \def\y@mathabcbox#1{\if :#1\else
   \y@vr\hbox to \y@b@xdim{\hfil$#1$\hfil}\fi}
\else
 \def\y@abcbox#1{\y@vr\hbox to \y@b@xdim{\hfil#1\hfil}}
 \def\y@mathabcbox#1{\y@vr\hbox to \y@b@xdim{\hfil$#1$\hfil}}
\fi
\def\y@setdim{%
  \ify@autoscale%
   \ifvoid1\else\typeout{Package youngtab: box1 not free! Expect an
     error!}\fi%
   \setbox1=\hbox{A}\y@b@xdim=1.6\ht1 \setbox1=\hbox{}\box1%
  \else\y@b@xdim=\y@boxdim \advance\y@b@xdim by -2\y@linethick
  \fi}
\newcount\y@counter
\newif\ify@islastarg
\def\y@lastargtest#1,#2 {\if\space #2 \y@islastargtrue
  \else\y@islastargfalse\fi}
\def\y@emptyboxes#1{\y@counter=#1\loop\ifnum\y@counter>0
  \advance\y@counter by -1 \y@emptybox\repeat}
\def\y@nelineemptyboxes#1{%
  \vbox{%
    \hrule height\y@linethick%
    \hbox{\y@emptyboxes{#1}\y@vr}
    \hrule height\y@linethick}\vspace{-\y@linethick}}
\def\yng(#1){%
  \y@setdim%
  \hspace{\y@interspace}%
  \ifmmode\ify@vcenter\vcenter\fi\fi{%
  \y@lastargtest#1,
  \vbox{\offinterlineskip
    \ify@islastarg
     \y@nelineemptyboxes{#1}
    \else
     \y@ungempty(#1)
    \fi}}\hspace{\y@interspace}}
\def\y@ungempty(#1,#2){%
  \y@nelineemptyboxes{#1}
  \y@lastargtest#2,
  \ify@islastarg
   \y@nelineemptyboxes{#2}
  \else
   \y@ungempty(#2)
  \fi}
\def\y@nelettertest#1#2. {\if\space #2 \y@islastargtrue
  \else\y@islastargfalse\fi}
\def\y@abcboxes#1#2.{%
  \ify@stdtext\y@abcbox#1\else\y@mathabcbox#1\fi%
  \y@nelettertest #2.
  \ify@islastarg\unskip%
   \ify@stdtext\y@abcbox{#2}\else\y@mathabcbox{#2}\fi%
  \else\y@abcboxes#2.\fi}
\ify@enable@skew
 \newdimen\y@full@b@xdim
 \newcount\y@m@veright@cnt
 \def\y@get@m@veright@cnt#1#2.{%
   \if :#1 \advance\y@m@veright@cnt by 1\y@get@m@veright@cnt#2.\fi}
 \let\y@setdim@=\y@setdim
 \def\y@setdim{%
   \y@setdim@ \y@full@b@xdim=\y@b@xdim
   \advance\y@full@b@xdim by 1\y@linethick}
 \def\y@m@veright@ifskew#1{
   \y@m@veright@cnt=0 \y@get@m@veright@cnt#1.
   \moveright \y@m@veright@cnt\y@full@b@xdim}
\else
 \def\y@m@veright@ifskew#1{}
\fi
\def\y@nelineabcboxes#1{%
  \y@nelettertest #1.
  \ify@islastarg
   \y@m@veright@ifskew{#1}
    \vbox{
      \hrule height\y@linethick%
      \hbox{\ify@stdtext\y@abcbox#1\else\y@mathabcbox#1\fi\y@vr}
      \hrule height\y@linethick}\vspace{-\y@linethick}
  \else
   \y@m@veright@ifskew{#1}
    \vbox{
      \hrule height\y@linethick%
      \hbox{\y@abcboxes #1.\y@vr}%
      \hrule height\y@linethick}\vspace{-\y@linethick}
  \fi}
\def\young(#1){%
  \y@setdim%
  \hspace{\y@interspace}%
  \y@lastargtest#1,
  \ifmmode\ify@vcenter\vcenter\fi\fi{%
  \vbox{\offinterlineskip
    \ify@islastarg\y@nelineabcboxes{#1}%
    \else\y@ungabc(#1)%
    \fi}}\hspace{\y@interspace}}
\def\y@ungabc(#1,#2){%
  \y@nelineabcboxes{#1}%
  \y@lastargtest#2,
  \ify@islastarg\y@nelineabcboxes{#2}%
  \else\y@ungabc(#2)%
  \fi}
\makeatother

\numberwithin{equation}{section}

\newcommand{\bC}{\mathbb{C}}

\newcommand{\bZ}{\mathbb{Z}}

\newcommand{\mf}[1]{\mathfrak{#1}}

\def\sgn{\mathrm{sgn}}

\def\gl{\mathfrak{gl}}
\def\sl{\mathfrak{sl}}

\def\miniyoung(#1){\scalebox{0.6}{$\young(#1)$}}
\def\mminiyoung(#1){\scalebox{0.4}{$\young(#1)$}}

\def\imm{{\textstyle \mathrm{Imm}_\lambda}}

\newcommand{\dete}[1]{{\textstyle \det_{#1}}}

\def\b1{\boldsymbol{1}}

\def\({ \left( }
\def\){ \right)}
\def\[{ \left[ }
\def\]{ \right]}

\theoremstyle{plain}
\newtheorem{thm}{Theorem}[section]
\newtheorem{prop}[thm]{Proposition}
\newtheorem{lem}[thm]{Lemma}
\newtheorem{cor}[thm]{Corollary}
\theoremstyle{definition}

\newtheorem{example}{Example}[section]
\newtheorem{remark}{Remark}[section]
\theoremstyle{conjecture}

\theoremstyle{problem}

\newcommand{\textem}[1]{{\bfseries #1}}

\title{\bfseries Alpha-determinant cyclic modules of $\mathfrak{gl}_n(\mathbb{C})$}
\author{\textsc{Sho MATSUMOTO}\footnote{Research Fellow of the Japan Society 
for the Promotion of Science, partially supported by Grant-in-Aid 
for Scientific Research (C) No. 17006193.}
\ and \
\textsc{Masato WAKAYAMA}\footnote{Partially supported by Grant-in-Aid 
for Scientific Research (B) No. 15340012.
}}

\date{August 19, 2005}

\begin{document}

\maketitle

\begin{abstract}
The alpha-determinant unifies and interpolates the notion of the determinant and permanent.
We determine the irreducible decomposition of 
the cyclic module of $\gl_n(\bC)$ defined by the alpha-determinant.
The degeneracy of the irreducible decomposition is determined
by the content polynomial of a given partition.

\par\noindent
\textem{Key Words} : 
Schur module, Young symmetrizer, permanent, Frobenius notation, irreducible decomposition,
cyclic module.
\end{abstract}


%
\section{Introduction}
%

Let $X=(x_{ij})_{1 \le i,j \le n}$ be a matrix with commutative variables $x_{ij}$.
For a complex number $\alpha$, 
the $\alpha$-determinant of $X$ is defined by
$$
\dete{\alpha}(X) = 
\sum_{\sigma \in \mf{S}_n} \alpha^{n-\nu_n(\sigma)} \prod_{i=1}^n x_{i \sigma(i)},
$$
where $\mf{S}_n$ is the symmetric group of degree $n$ and 
$\nu_n(\sigma)$ stands for the number of cycles in the cycle decomposition of a permutation 
$\sigma \in \mf{S}_n$.
The $\alpha$-determinant is nothing but 
the permanent if $\alpha=1$ and the (usual) determinant if $\alpha=-1$,
and hence, it interpolates these two.
It appears as the coefficient in the Taylor expansion of the power $\det(I-\alpha X)^{-1/\alpha}$
of the characteristic polynomial of $X$
and defines a generalization of the boson, poisson and fermion point processes,
see \cite{ST, V}.
Also, its Pfaffian analogue has been developed in \cite{AlphaPf}.

It is a natural question whether the $\alpha$-determinant can be interpreted
as an invariant
like the usual determinant
(and also the $q$-determinant in quantum group theory).
Denote by $\mathcal{P}(\mathrm{Mat}_{n \times n})$ the ring  of polynomials 
in variables $\{x_{ij}\}_{1 \le i,j \le n}$.
Let $\{ E_{ij} \}_{1 \le i,j \le n}$ be the natural basis of the Lie algebra $\mf{g}=\gl_n(\bC)$.
When $n=2$ and $\alpha \not=0$,
consider the linear map $\rho^{(\alpha)}_2$ of $\mf{gl}_2(\bC)$ on $\mathcal{P}(\mathrm{Mat}_{n \times n})$
determined by
\begin{align*}
\rho^{(\alpha)}_2(E_{11}) =& x_{11}\partial_{11} + x_{12}\partial_{12}, &
\rho^{(\alpha)}_2(E_{12}) =& \frac{1}{\sqrt{-\alpha}} 
\(x_{11} \partial_{21} -\alpha x_{12}\partial_{22} \), \\
\rho^{(\alpha)}_2(E_{21}) =& \frac{1}{\sqrt{-\alpha}}
\(-\alpha x_{21} \partial_{11} + x_{22} \partial_{12} \), &
\rho^{(\alpha)}_2(E_{22}) =& x_{21}\partial_{21} + x_{22}\partial_{22},
\end{align*}
where $\partial_{ij}=\frac{\partial}{\partial x_{ij}}$.
Then $\rho^{(\alpha)}_2$ defines a representation of $\mf{gl}_2(\bC)$ and 
$\rho^{(\alpha)}_2 (E_{ii}) \dete{\alpha}(X) =\dete{\alpha}(X)$, 
$\rho^{(\alpha)}_2 (E_{ij}) \dete{\alpha}(X) =0$ for $i \not= j$.
This is not, however, true for $n \ge 3$.
Precisely, 
although the map $\rho^{(\alpha)}_n$ given by 
$\rho^{(\alpha)}_n(E_{ij})= \sum_{k=1}^n \beta^{|i-k| - |j-k|} x_{ik} \partial_{jk}$, 
where $\beta =\sqrt{-\alpha}$, defines a representation of $\mf{gl}_n(\bC)$ 
on $\mathcal{P}(\mathrm{Mat}_{n \times n})$,
$\rho_n^{(\alpha)}(E_{ij}) \dete{\alpha}(X) \not= 0 \ (i \not= j)$ in general.
(One can actually show that $\rho_n^{(\alpha)}$ is equivalent to 
the usual action of $\mf{gl}_n(\bC)$ determined by 
$\rho(E_{ij}) =\rho^{(-1)}_n(E_{ij}) =\sum_{k=1}^n x_{ik} \partial_{jk}$ 
on $\mathcal{P}(\mathrm{Mat}_{n \times n})$ when $\alpha \not=0$.
Indeed, the map $f(x_{ij}) \mapsto f(\beta^{|i-j|} x_{ij})$ 
is the intertwining operator from $(\rho, \mathcal{P}(\mathrm{Mat}_{n \times n}))$ 
to $(\rho_n^{(\alpha)}, \mathcal{P}(\mathrm{Mat}_{n \times n}))$.)

Then, a question which subsequently arises
is what the structure of the smallest invariant subspace of $\mathcal{P}(\mathrm{Mat}_{n \times n})$
which contains $\dete{\alpha}(X)$ is.
Thus, the aim of the present paper is to investigate
a cyclic module $V_n^{(\alpha)} =U(\mf{g}) \dete{\alpha}(X)$,
under the representation
$\rho$ on $\mathcal{P}(\mathrm{Mat}_{n \times n})$
for each $\alpha \in \bC$.
Clearly, $V_n^{(-1)}$ is the one-dimensional determinant representation.

We adopt the notations for partitions used in \cite{Mac} 
and for representations of $\mf{gl}_n (\bC)$ used in \cite{Fulton} and \cite{Weyl}.
A partition $\lambda$ is a weakly decreasing sequence 
$\lambda=(\lambda_1, \lambda_2, \dots)$
of non-negative integers such that $\lambda_j =0$ for sufficiently large $j$.
We usually identify a partition $\lambda$ with the corresponding Young diagram.
Write $\lambda \vdash n$ if $\sum_{j \ge 1} \lambda_j=n$
and denote by $\lambda'=(\lambda_1', \lambda_2', \dots)$ the conjugate partition of $\lambda$. 
Let $E^\lambda$ denote the Schur module (or called the Weyl module)
corresponding to $\lambda$
and $f^{\lambda}$ the number of standard tableaux of shape $\lambda$.

The following is our main result,
which describes the irreducible decomposition of $V_n^{(\alpha)}$.

\begin{thm} \label{MainTheorem}
For $k=1,2, \dots, n-1$,
\begin{equation}
V_n^{(\frac{1}{k})} \cong 
\bigoplus_{\begin{subarray}{c} \lambda \vdash n, \\ \lambda'_1 \leq k \end{subarray}} 
(E^{\lambda})^{\oplus f^{\lambda}}
\qquad \text{and} \qquad
V_n^{(-\frac{1}{k})} \cong
\bigoplus_{\begin{subarray}{c} \lambda \vdash n, \\ \lambda_1 \leq k \end{subarray}} 
(E^{\lambda})^{\oplus f^{\lambda}}.
\end{equation}
For $\alpha \in \bC \setminus \{\pm 1, \pm \frac{1}{2}, \dots, \pm \frac{1}{n-1} \}$,
\begin{equation} \label{EqDecompositionGeneric}
V_n^{(\alpha)} \cong (\bC^n)^{\otimes n} \cong
\bigoplus_{\lambda \vdash n}
(E^{\lambda})^{\oplus f^{\lambda}}.
\end{equation}

\end{thm}

\begin{example}
When $n=3$ the irreducible decomposition of $V_n^{(\alpha)}$ is given by
$$
V_3^{(\alpha)} \cong
\begin{cases} 
E^{(3)} & \text{if $\alpha=1$}, \\
E^{(3)} \oplus E^{(2,1)} \oplus E^{(2,1)} & \text{if $\alpha=\frac{1}{2}$}, \\
E^{(1,1,1)} & \text{if $\alpha=-1$}, \\
E^{(2,1)} \oplus E^{(2,1)} \oplus E^{(1,1,1)}  & \text{if $\alpha=-\frac{1}{2}$}, \\
E^{(3)} \oplus E^{(2,1)} \oplus E^{(2,1)} \oplus E^{(1,1,1)} & \text{otherwise}. \qed
\end{cases}
$$
\end{example}

We remark that each Schur module possesses a canonial basis formed by $\alpha$-determinants
(Theorem 3.3).
Although the $\alpha$-determinant is not an invariant,
this fact implies that it has a rich symmetry.

%
\section{The $U(\mf{gl}_n)$-module $V_n^{(\alpha)}$}
%

Let $V_n^{(\alpha)}= \rho(U(\mf{gl}_n)) \dete{\alpha}(X)$.
Put $[n]=\{1, 2, \dots, n\}$ and
$$
D^{(\alpha)}(i_1, i_2, \dots, i_n)
=\dete{\alpha} 
\begin{pmatrix} 
x_{i_1 1} & x_{i_1 2} & \ldots & x_{i_1 n} \\
x_{i_2 1} & x_{i_2 2} & \ldots & x_{i_2 n} \\
\vdots & \vdots & \ddots & \vdots \\
x_{i_n 1} & x_{i_n 2} & \ldots & x_{i_n n} 
\end{pmatrix}
$$ 
for any $i_1, \dots, i_n \in [n]$. 
In particular, $\dete{\alpha}(X)=D^{(\alpha)}(1,2, \dots, n)$.
We abbreviate $\rho(E_{ij})$ to $E_{ij}$ for simplicity.
When there is no fear of confusion, we abbreviate $\nu_n$ to $\nu$.

\begin{lem} \label{LemACTION}
\begin{equation}
E_{pq} \cdot D^{(\alpha)}(i_1, \dots, i_n) = 
\sum_{k=1}^n \delta_{i_k, q} D^{(\alpha)}(i_1, \dots, i_{k-1}, p, i_{k+1}, \dots, i_n).
\end{equation}
\end{lem}

\begin{proof}
It is straightforward.
In fact,
\begin{align*}
E_{pq} \cdot D^{(\alpha)}(i_1, \dots, i_n) 
=& 
\sum_{j=1}^n x_{pj} \frac{\partial}{\partial x_{qj}}
\sum_{\sigma \in \mf{S}_n} \alpha^{n-\nu(\sigma)} 
x_{i_1 \sigma(1)} \cdots x_{i_n \sigma(n)} \\
=&
\sum_{j=1}^n  \sum_{\sigma \in \mf{S}_n} \alpha^{n-\nu(\sigma)} 
\sum_{k=1}^n x_{pj} \delta_{i_k,q} \delta_{\sigma(k),j} 
x_{i_1 \sigma(1)} \cdots \widehat{x_{i_k \sigma(k)}} \cdots x_{i_n \sigma(n)} \\
=&
\sum_{k=1}^n  \delta_{i_k,q} \sum_{\sigma \in \mf{S}_n} \alpha^{n-\nu(\sigma)} 
x_{p \sigma(k)}
x_{i_1 \sigma(1)} \cdots \widehat{x_{i_k \sigma(k)}} \cdots x_{i_n \sigma(n)} \\
=& \sum_{k=1}^n \delta_{i_k, q} D^{(\alpha)}(i_1, \dots, i_{k-1}, p, i_{k+1}, \dots, i_n),
\end{align*}
where $\widehat{x_{kl}}$ stands for the omission of $x_{kl}$. 
\end{proof}

\begin{example}
We see that
$E_{21} \cdot D^{(\alpha)} (4,1,2,1) = D^{(\alpha)} (4,2,2,1) +  D^{(\alpha)} (4,1,2,2)$,
$E_{11} \cdot D^{(\alpha)} (4,1,2,1) = 2 D^{(\alpha)} (4,1,2,1)$, 
and $E_{43} \cdot D^{(\alpha)} (4,1,2,1) =0$. \qed
\end{example}

The symmetric group $\mf{S}_n$ acts also on $V_n^{(\alpha)}$ from the right
by
$D^{(\alpha)}(i_1, \dots, i_n) \cdot \sigma 
= D^{(\alpha)} (i_{\sigma(1)}, \dots, i_{\sigma(n)})$.

\begin{lem}  \label{LemSPAN}
The space $V_n^{(\alpha)}$ is the complex vector space spanned by
$\{D^{(\alpha)}(i_1, \dots, i_n) \ | \ i_1, \dots, i_n \in [n]\}$.
\end{lem}

\begin{proof}
Since the vector space spanned by all $D^{(\alpha)}(i_1, \dots, i_n)$  
contains $V_n^{(\alpha)}$ by Lemma \ref{LemACTION},
we prove that all $D^{(\alpha)}(i_1, \dots, i_n)$ are contained in $V_n^{(\alpha)}$.
For $1 \le p < q \le n$, we have
$$
V_n^{(\alpha)} \ni (E_{pq} E_{qp}-1) \cdot D^{(\alpha)} (1,2, \dots, n) = 
D^{(\alpha)} (\tau(1), \dots, \tau(n)),
$$
where $\tau$ is the transposition $(p,q)$ of $p$ and $q$.
It follows that, for each $\sigma \in \mf{S}_n$,  
$D^{(\alpha)}(\sigma(1), \dots, \sigma(n))=Y \cdot D^{(\alpha)}(1,\dots,n)$
for some $Y=Y_\sigma \in U(\mf{g})$.
For any $1 \le i_1 \le i_2 \le \dots \le i_n \le n$,
suppose there exists $X=X_{i_1 \dots i_n} \in U(\mf{g})$ such that
$D^{(\alpha)}(i_1, \dots, i_n)=X \cdot D^{(\alpha)}(1, \dots, n)$.
For any $j_1, \dots, j_n \in [n]$,
we have $D^{(\alpha)}(j_1, \dots, j_n) = D^{(\alpha)}(i_{\sigma(1)}, \dots, i_{\sigma(n)})$  
for some $\sigma \in \mf{S}_n$ and $i_1 \le \dots \le i_n$.
Hence, since the action of $\mf{gl}_n(\bC)$ and of $\mf{S}_n$ commute,
we see that
\begin{align*}
& D^{(\alpha)}(j_1, \dots, j_n)= D^{(\alpha)}(i_1, \dots, i_n) \cdot \sigma \\
=& \( X \cdot D^{(\alpha)}(1, \dots, n)\) \cdot \sigma
= X \cdot \( D^{(\alpha)}(1, \dots, n) \cdot \sigma \) 
= X Y \cdot D^{(\alpha)}(1, \dots, n)
\end{align*}
and $D^{(\alpha)}(j_1, \dots, j_n)$ is contained in $V_n^{(\alpha)}$.
Therefore it is sufficient to prove $D^{(\alpha)} (i_	1, \dots, i_n) \in V_n^{(\alpha)}$
for $i_1 \le \dots \le i_n$.

For any sequence $(i_1, \dots, i_n)$ such that $i_k \le k$ for any $k$,
we have
$D^{(\alpha)} (i_1, i_2, \dots, i_n) = E_{i_n n} \cdots E_{i_2 2} E_{i_1 1} 
\cdot D^{(\alpha)}(1,2,\dots, n) \in V_n^{(\alpha)}$.
In fact, by Lemma \ref{LemACTION}, we see that
$E_{i_k k} \cdot D^{(\alpha)} ( i_1, \dots, i_{k-1}, k, k+1, \dots, n)
= D^{(\alpha)} ( i_1, \dots, i_{k-1}, i_k, k+1, \dots, n)$
because $i_1 \le \dots \le i_{k-1} < k$ for any $1 \le k \le n$.
Suppose there exists $k$ such that $i_j \le j$ for any $j < k$ and $i_k >k$.
We prove $D^{(\alpha)} (i_1, \dots, i_n) \in V_n^{(\alpha)}$
for such sequences $(i_1, \dots, i_n)$
by induction with respect to the lexicographic order.
Since $i_1 \le \dots \le i_{k-1} < k \le i_k-1 < i_{k+1} \le \dots \le i_n$
and
$D^{(\alpha)} (i_1, \dots, i_{k-1}, i_k-1, i_{k+1}, \dots, i_n) \in V_n^{(\alpha)}$
by the induction assumption,
we have
$D^{(\alpha)}(i_1, \dots, i_{k}, \dots, i_n) 
= E_{i_k, i_k-1} \cdot D^{(\alpha)}(i_1, \dots, i_k-1, \dots, i_n) \in V_n^{(\alpha)}$
by Lemma \ref{LemACTION}.
Hence we obtain our claim.
\end{proof}

The universal enveloping algebra $U(\mf{g})$ of $\mf{g}=\mf{gl}_n(\bC)$ acts on the $n$-tensor product  
$(\bC^n)^{\otimes n}=\bC^n \otimes \cdots \otimes \bC^n$ from the left  by
$$
E_{pq} \cdot (\boldsymbol{e}_{i_1} \otimes  \cdots \otimes \boldsymbol{e}_{i_n})
= \sum_{k=1}^n \boldsymbol{e}_{i_1} \otimes \cdots \otimes E_{pq} \boldsymbol{e}_{i_k} 
\otimes \cdots \otimes \boldsymbol{e}_{i_n}
= \sum_{k=1}^n \delta_{i_k, q} \boldsymbol{e}_{i_1} \otimes \cdots \otimes  \boldsymbol{e}_{p} 
\otimes \cdots \otimes \boldsymbol{e}_{i_n},
$$
where $\{\boldsymbol{e}_k\}_{k=1}^n$ is the natural basis of $\bC^n$.
From this fact together with Lemma \ref{LemACTION} and Lemma \ref{LemSPAN},
we have the

\begin{prop} \label{PropProjection}
Let $\Phi^{(\alpha)}_n$ be the linear map from $(\bC^n)^{\otimes n}$
to $V_n^{(\alpha)}$ defined by
$$
\Phi^{(\alpha)}_n (\boldsymbol{e}_{i_1} \otimes  \cdots \otimes \boldsymbol{e}_{i_n})
= D^{(\alpha)} (i_1, \dots, i_n)
$$
for each $i_1, \dots, i_n \in [n]$.
Then $\Phi^{(\alpha)}_n$ is a $U(\mf{g})$-module homomorphism.
In particular, $V_n^{(\alpha)}$ is isomorphic to 
a quotient module $(\bC^n)^{\otimes n} \big/ \mathrm{Ker\,} \Phi_n^{(\alpha)}$ of $(\bC^n)^{\otimes n}$.
\qed
\end{prop}

Notice that, when $\alpha=0$,
the homomorphism 
$\Phi^{(0)}_n (\boldsymbol{e}_{i_1} \otimes  \cdots \otimes \boldsymbol{e}_{i_n})
= D^{(0)} (i_1, \dots, i_n)= x_{i_1 1} \cdots x_{i_n n}$
is clearly bijective,
and therefore
$V_n^{(0)} \cong (\bC^n)^{\otimes n}$.
Hence
we have the irreducible decomposition of $V_n^{(0)}$ as
\begin{equation} \label{EqDecompositionZero}
V_n^{(0)} \cong \bigoplus_{\lambda \vdash n} (E^{\lambda})^{\oplus f^{\lambda}}.
\end{equation}
The symmetric group $\mf{S}_n$ acts on $(\bC^n)^{\otimes n}$ from the right by
$(\boldsymbol{e}_{i_1} \otimes  \cdots \otimes \boldsymbol{e}_{i_n}) \cdot \sigma
= \boldsymbol{e}_{i_{\sigma(1)}} \otimes  \cdots \otimes \boldsymbol{e}_{i_{\sigma(n)}}$
for any $\sigma \in \mf{S}_n$.

%
\section{A formula for the number of cycles}
%

A numbering of shape $\lambda \vdash n$ is a way of putting distinct elements in $[n]$
in each box of the Young diagram $\lambda$.
Let $R(T)$ be the row group (or called the Young subgroup) 
of a numbering $T$, i.e., permutations in $R(T)$
permutate the entries of each row among themselves.
The column group $C(T)$ is also defined similarly.

Recall the Frobenius notation $(a_1, a_2, \dots, a_d | b_1, b_2, \dots, b_d)$
of a partition $\lambda$, where
$a_i= \lambda_i-i \ge 0$ and $b_i=\lambda_i'-i \ge 0$ for $1 \le i \le d$.
Then the content polynomial $f_\lambda(\alpha)$ (\cite{Mac}) for the partition $\lambda$
is written as
\begin{equation}
f_{\lambda} (\alpha)=
\prod_{i=1}^d \left\{ \prod_{j=1}^{a_i}(1+ j \alpha) \cdot 
\prod_{j=1}^{b_i} (1-j\alpha) \right\}.
\end{equation}
Note that $f_{\lambda}(\alpha)$ satisfies 
$f_{\lambda}(\alpha)=f_{\lambda'}(-\alpha)$.
We have the following formula for the number $\nu=\nu_n$ of cycles.

\begin{prop} \label{ThmCycleFormula}
Let $T$ be a numbering of shape $\lambda \vdash n$.
Then
\begin{align}
 & \sum_{q \in C(T)} \sgn(q) \sum_{p \in R(T)} \alpha^{n-\nu(p q \sigma)} \label{EqCycleFormula} \\
=& \begin{cases} \sgn(q_0) f_{\lambda}(\alpha) 
  & \text{if $\sigma=q_0 p_0$ for some $q_0 \in C(T)$ and $p_0 \in R(T)$}, \\
  0 & \text{otherwise}. 
\end{cases} \notag
\end{align}
\end{prop}

\begin{example} 
When $T=\miniyoung(12,3)$ we have
$$
\sum_{q \in C(T)} \sgn(q) \sum_{p \in R(T)} \alpha^{3-\nu(pq \sigma)}
= \begin{cases}
(1+\alpha)(1-\alpha) & \text{for $\sigma=(1)$ or $(12)$}, \\
-(1+\alpha)(1-\alpha) & \text{for $\sigma=(13)$ or $(123)$}, \\
0 & \text{for $\sigma=(23)$ or $(132)$}. \qed
\end{cases} \
$$
\end{example}

\begin{example}
For
$T=\young({1}{2}{\cdots}{n})$
Proposition \ref{ThmCycleFormula} says
\begin{equation} \label{EqStanley}
\sum_{\sigma \in \mf{S}_n} \alpha^{n-\nu(\sigma)} = \prod_{j=1}^{n-1} (1+ j \alpha). \qed
\end{equation}
\end{example}

Though we can give an even simpler proof of this proposition 
based on a well-known representation theory of $\mf{S}_n$ and $GL_m$ 
(see, Section \ref{SectionQuantum} and  \ref{SectionImmanant}),
since this does not work in the quantum groups situation,
we provide here a direct and elementary proof.

\begin{proof}[Proof of Proposition \ref{ThmCycleFormula}.]
First we prove the assertion when $T$ is a standard tableau.
Suppose that the $(t,s)$-box in $T$ is numbered by $n$.
Let $\{x_1, \dots, x_{s-1},n \}$ and $\{y_1, \dots, y_{t-1},n \}$ be
the entries on the $t$-th row and $s$-th column in $T$, respectively.
Let $\tilde{T}$ be the standard tableau obtained by removing the $(t,s)$-box from $T$
and let $\tilde{\lambda}$ be the shape of $\tilde{T}$.
For any $\sigma \in \mf{S}_n$, 
we define $\tilde{\sigma}$ by
$\tilde{\sigma}=\sigma \cdot (\sigma^{-1}(n),n)$.
The map $\sigma \mapsto \tilde{\sigma}$ is a projection 
from $\mf{S}_{n}$ to $\mf{S}_{n-1}$.
Then 
\begin{align}
& \sum_{q \in C(T)} \sgn (q) \sum_{p \in R(T)} 
\alpha^{n-\nu_n(pq \sigma)} \label{EqEReduction}  \\
=& F(\sigma;n,n)+ \sum_{k=1}^{s-1} F(\sigma; x_k,n) + \sum_{l=1}^{t-1} F(\sigma; n, y_l)
+\sum_{k=1}^{s-1} \sum_{l=1}^{t-1} F(\sigma; x_k,y_l), \notag
\end{align}
where
\begin{align*}
F(\sigma;n,n) =& \sum_{q \in C(\tilde{T}) (n)} \sgn (q) \sum_{p \in (n) R(\tilde{T})} 
\alpha^{n-\nu_n(pq \sigma)} 
= \sum_{\tilde{q} \in C(\tilde{T})} \sgn (\tilde{q}) \sum_{\tilde{p} \in R(\tilde{T})} 
\alpha^{n-\nu_n(\tilde{p} \tilde{q} \sigma )}, \\
F(\sigma; x_k,n) =& \sum_{q \in C(\tilde{T}) (n)} \sgn (q) \sum_{p \in (x_k,n) R(\tilde{T})} 
\alpha^{n-\nu_n(pq \sigma)} 
= \sum_{\tilde{q} \in C(\tilde{T})} 
  \sgn (\tilde{q}) \sum_{\tilde{p} \in R(\tilde{T})} 
  \alpha^{n-\nu_n((x_k,n) \tilde{p} \tilde{q} \sigma)}, \\
F(\sigma; n, y_l) =& \sum_{q \in C(\tilde{T}) (y_l,n)} \sgn (q) \sum_{p \in (n) R(\tilde{T})} 
  \alpha^{n-\nu_n(pq \sigma)}
= - \sum_{\tilde{q} \in C(\tilde{T})} \sgn (\tilde{q}) 
   \sum_{\tilde{p} \in R(\tilde{T})} 
   \alpha^{n-\nu_n(\tilde{p} \tilde{q}(y_l,n) \sigma)}, \\
F(\sigma; x_k, y_l) =& \sum_{q \in C(\tilde{T}) (y_l,n)} \sgn (q) \sum_{p \in (x_k,n) R(\tilde{T})} 
  \alpha^{n-\nu_n(pq \sigma)}
= -\sum_{\tilde{q} \in C(\tilde{T})} \sgn (\tilde{q}) \sum_{\tilde{p} \in R(\tilde{T})} 
  \alpha^{n-\nu_n((x_k, n) \tilde{p}\tilde{q} (y_l,n) \sigma)}.
\end{align*}

We prove \eqref{EqCycleFormula} by induction on $n$
for five cases:
(i) $\sigma$ is the identity; 
(ii) $\sigma=q_0 p_0$ for some $q_0 \in C(T)$ and $p_0 \in R(T)$;
(iii) $\sigma=(z,n)$ for $z \in [n-1] \setminus \{x_1, \dots, x_{s-1}, y_1, \dots, y_{t-1}\} $;
(iv)  $\sigma$ is not of the form $\sigma=q_0 p_0$ for any $q_0 \in C(T)$ and $p_0 \in R(T)$
but $\tilde{\sigma}=\tilde{q}_0 \tilde{p}_0$ 
for some $\tilde{q}_0 \in C(\tilde{T})$ and $\tilde{p}_0 \in R(\tilde{T})$;
(v) $\sigma$ is not of the form $\sigma=q_0 p_0$ for any $q_0 \in C(T)$ and $p_0 \in R(T)$,
and $\tilde{\sigma}$ is not of the form 
$\tilde{\sigma}=\tilde{q}_0 \tilde{p}_0$ 
for any $\tilde{q}_0 \in C(\tilde{T})$ and $\tilde{p}_0 \in R(\tilde{T})$.
Note that (i) (resp. (iii)) is the special case in (ii) (resp. (iv)).
Assume that \eqref{EqCycleFormula} holds for $\tilde{T}$ and $\tilde{\sigma}$.

Case (i).
We have 
$$
F((1);n,n) =\sum_{\tilde{q} \in C(\tilde{T})} \sgn (\tilde{q}) \sum_{\tilde{p} \in R(\tilde{T})} 
\alpha^{n-1-\nu_{n-1}(\tilde{p}  \tilde{q})} =f_{\tilde{\lambda}}(\alpha)
$$
by the assumption of the induction for $\tilde{\sigma}=(1)$
because $\nu_n(\tilde{p}\tilde{q} ) = \nu_{n-1}( \tilde{p}\tilde{q} )+1$.
We have also
$$
F((1);x_k,n) 
= \alpha \sum_{\tilde{q} \in C(\tilde{T})} \sgn (\tilde{q}) \sum_{\tilde{p} \in R(\tilde{T})} 
\alpha^{n-1-\nu_{n-1}(\tilde{p}  \tilde{q})} = \alpha f_{\tilde{\lambda}}(\alpha)
$$
by the induction hypothesis again 
because $\nu_n((x_k,n) \tilde{p} \tilde{q}) =\nu_{n-1}(\tilde{p} \tilde{q})$.
Similarly, $F((1); n,y_l)= - \alpha f_{\tilde{\lambda}}(\alpha)$.
Further,
$$
F((1); x_k,y_l)
= \sum_{\tilde{q} \in C(\tilde{T})} \sgn (\tilde{q}) \sum_{\tilde{p} \in R(\tilde{T})} 
\alpha^{n-\nu_n((x_k,y_l, n)  \tilde{p}\tilde{q})}  \\
= \alpha \sum_{\tilde{q} \in C(\tilde{T})} \sgn (\tilde{q}) \sum_{\tilde{p} \in R(\tilde{T})} 
\alpha^{n-1-\nu_{n-1}((x_k,y_l) \tilde{p}  \tilde{q})}
$$
and this equals zero because
permutation $(x_k, y_l)$ is not given in the form of
a product of an element of $C(\tilde{T})$ and one of $R(\tilde{T})$.
Hence by \eqref{EqEReduction} we obtain
\begin{equation*} 
\sum_{q \in C(T)} \sgn (q) \sum_{p \in R(T)} 
\alpha^{n-\nu_n(pq)} = (1+(s-t) \alpha)f_{\tilde{\lambda}}(\alpha).
\end{equation*}
The Frobenius notation of $\tilde{\lambda}$ is given by
$(a_1, \dots, a_{d-1} | b_1, \dots, b_{d-1})$ if $s=t$; or
$(a_1, \dots, a_t-1, \dots, a_d | b_1, \dots, b_d)$ and $a_t=s-t$ if $s >t$; or
$(a_1, \dots, a_d | b_1, \dots, b_s-1, \dots b_d)$ and $b_s=t-s$ if $s<t$.
Therefore $f_{\lambda}(\alpha)=(1+(s-t) \alpha)f_{\tilde{\lambda}}(\alpha)$,
and \eqref{EqCycleFormula} holds for the case (i).

Case (ii).
Since
$$
\sum_{q \in C(T)} \sgn (q) \sum_{p \in R(T)} 
\alpha^{n-\nu_n(pq q_0 p_0)} = \sgn(q_0) \sum_{q \in C(T)} \sgn (q) \sum_{p \in R(T)} 
\alpha^{n-\nu_n(pq)},
$$
the equality \eqref{EqCycleFormula} holds for the case (ii) by the result for (i).

Case (iii).
Let $\sigma=(z,n)$ and suppose that the $(t',s')$-box in $T$ is numbered by $z$.
Since $T$ is a standard tableau, it holds $t'<t$ or $s' <s$.
We have $F(\sigma; n,n)= \alpha f_{\tilde{\lambda}}(\alpha)$
and
$$
F(\sigma;x_k,n)= \sum_{\tilde{q} \in C(\tilde{T})} 
  \sgn (\tilde{q}) \sum_{\tilde{p} \in R(\tilde{T})} 
  \alpha^{n-\nu_n(\tilde{p} \tilde{q} (x_k, z, n)) }
=\alpha \sum_{\tilde{q} \in C(\tilde{T})} 
  \sgn (\tilde{q}) \sum_{\tilde{p} \in R(\tilde{T})} 
  \alpha^{n-1-\nu_{n-1}( \tilde{p} \tilde{q} (x_k, z))}.
$$
If $s'<s$ then $(x_{s'},z)$ belongs to $C(\tilde{T})$
and $(x_k,z)$ is not given as a product of one in $C(\tilde{T})$ and one in $R(\tilde{T})$
for $k \not= s'$.
Therefore it follows by the assumption of the induction that
$$
\sum_{k=1}^{s-1} F(\sigma;x_k,n) =
\begin{cases} 
-\alpha f_{\tilde{\lambda}}(\alpha) & \text{if $s' <s$}, \\
0 & \text{otherwise}.
\end{cases}
$$
Similarly, we have
$$
\sum_{l=1}^{t-1} F(\sigma;n,y_l) =
\begin{cases} 
-\alpha f_{\tilde{\lambda}}(\alpha) & \text{if $t' <t$}, \\
0 & \text{otherwise}
\end{cases}
$$
and 
$$
\sum_{k=1}^{s-1} \sum_{l=1}^{t-1} F(\sigma;x_k,y_l) =
\begin{cases} 
\alpha f_{\tilde{\lambda}}(\alpha) & \text{if $s' < s$ and $t' <t$}, \\
0 & \text{otherwise}.
\end{cases}
$$
Hence, by \eqref{EqEReduction} we obtain 
$$
\sum_{q \in C(T)} \sgn (q) \sum_{p \in R(T)} 
\alpha^{n-\nu_n(pq (z,n))} =0. 
$$
This shows \eqref{EqCycleFormula} for (iii).

Case (iv).
In this case, $\sigma$ is expressed as 
$\sigma=(z',n) \tilde{q}_0 \tilde{p}_0= \tilde{q}_0 \tilde{p}_0 (z,n)$,
where $z' \in [n-1] \setminus \{y_1, \dots, y_{t-1} \}$ 
and $z \in [n-1] \setminus \{x_1, \dots, x_{s-1} \}$. 
Since
$$
\sum_{q \in C(T)} \sgn (q) \sum_{p \in R(T)} \alpha^{n- \nu_n(p q \tilde{q}_0 \tilde{p}_0 (z,n))}
=\sgn (q_0) \sum_{q \in C(T)} \sgn (q) \sum_{p \in R(T)} \alpha^{n- \nu_n(p q \tilde{p}_0 (z,n))},
$$
we may assume $\tilde{q}_0$ is the identity without loss of generality.
By a similar discussion, we may further assume $\tilde{p}_0$ is the identity,
and hence the case (iv) is reduced to (iii).

Case (v).
It is clear that $F(\sigma; n,n)$ is equal to just the following or the $\alpha$ times of
$$
\sum_{\tilde{q} \in C(\tilde{T})} \sgn (\tilde{q}) \sum_{\tilde{p} \in R(\tilde{T})} 
\alpha^{n-1-\nu_{n-1}(\tilde{p}  \tilde{q} \tilde{\sigma})}
$$
but this equals zero by the induction assumption.
Since $\nu_n((x_k,n) \tilde{q}_0 \tilde{p}_0 \tilde{\sigma}) 
= \nu_{n}(\tilde{q}_0 \tilde{p}_0 \tilde{\sigma}) \pm 1$
we obtain
$F(\sigma; x_k,n)=0$, 
and similarly  $F(\sigma;n,y_l)=F(\sigma; x_k,y_l)=0$.
Therefore \eqref{EqCycleFormula} holds for (v).
This completes the proof of the theorem for a standard tableau $T$.

In general, a numbering $T$ is written as $T= \tau \cdot U$
for some $\tau \in \mf{S}_n$ and a standard tableau $U$.
Here $\tau \cdot U$ is the numbering obtained from replacing $i$ by $\tau(i)$ 
in each box of $U$.
Since $R(T)= \tau R(U) \tau^{-1}$ and $C(T)= \tau C(U) \tau^{-1}$
it is proved \eqref{EqCycleFormula} for a numbering $T$ by the result for a standard tableau $U$.
\end{proof}

For a sequence $(i_1,\dots, i_n) \in [n]^n$ and a numbering $T$,
we define
\begin{equation} \label{EqDefBASIS}
v_{T}^{(\alpha)} (i_1, \dots, i_n) =D^{(\alpha)}(i_1, \dots, i_n) \cdot c_{T}
= \sum_{q \in C(T)} \sgn (q)  \sum_{p \in R(T)} 
D^{(\alpha)}(i_{qp(1)}, \dots, i_{qp(n)}),
\end{equation}
where $c_T= \sum_{q \in C(T)} \sgn (q)  \sum_{p \in R(T)}  q p$ 
is the Young symmetrizer.

\begin{cor} \label{ThmFactorization}
Let $\lambda \vdash n$.
For each sequence $(i_1, \dots, i_n) \in [n]^n$ and numbering $T$ of shape $\lambda$,
we have
\begin{equation} \label{EqFactorization}
v_{T}^{(\alpha)} (i_1, \dots, i_n)
= f_{\lambda}(\alpha)
v_{T}^{(0)} (i_1, \dots, i_n).
\end{equation}
\end{cor}

\begin{proof}
For a numbering $T$ of shape $\lambda$,
denote by $W^{(\alpha)}_T$ the space  
spanned by $\{ v_{T}^{(\alpha)} (i_1, \dots, i_n) \ | \ i_1, \dots, i_n \in [n] \}$. 
Since the Schur module $E^{\lambda}$ is isomorphic to the image 
of the map $(\bC^n)^{\otimes n} \rightarrow (\bC^n)^{\otimes n}$ given by $c_T$
for any numbering $T$, 
it follows from Proposition \ref{PropProjection} and \eqref{EqDefBASIS} that
$W_T^{(\alpha)}$ is isomorphic to  $\{0\}$ or $E^{\lambda}$.
If \eqref{EqFactorization} is proved  for a certain sequence $(j_1, \dots, j_n)$
satisfying $v_{T}^{(0)}(j_1, \dots, j_n) \not=0$, 
\eqref{EqFactorization} holds for any $(i_1, \dots, i_n) \in [n]^n$
because $W_T^{(\alpha)}$ is the cyclic module $U(\mf{g}) v_{T}^{(\alpha)}(j_1, \dots, j_n)$ 
and the action of $\mf{g}$ is independent of $\alpha$.

We prove \eqref{EqFactorization}
for the case where $(i_1, i_2, \dots, i_n)= (1,2,\dots, n)$.
Then we have
\begin{align*}
v_{T}^{(\alpha)} (1,\dots,n) 
=& \sum_{q \in C(T)} \sgn(q) \sum_{p \in R(T)} \sum_{\sigma \in \mf{S}_n}
\alpha^{n-\nu(\sigma)} x_{q p \sigma(1),1} \cdots x_{q p \sigma(n),n} \\
=&  \sum_{\sigma \in \mf{S}_n} \( \sum_{q \in C(T)} \sgn(q) \sum_{p \in R(T)}
\alpha^{n-\nu(p q \sigma)}  \)x_{\sigma(1),1} \cdots x_{\sigma(n),n}.
\end{align*}
By Proposition \ref{ThmCycleFormula} we see that
$$
v_{T}^{(\alpha)}(1,\dots,n) = f_{\lambda}(\alpha)
\sum_{q_0 \in C(T)} \sgn(q_0) \sum_{p_0 \in R(T)}
x_{q_0 p_0(1),1} \cdots x_{q_0 p_0(n),n}
=f_{\lambda}(\alpha) v_{T}^{(0)} (1,\dots,n).
$$
If $q_0 \not= q_0'$ or $p_0 \not= p'_0$ then 
$q_0 p_0 \not= q_0' p_0'$.
Indeed, if $q_0 p_0 = q_0' p_0'$ then $C(T) \ni (q_0')^{-1} q_0 = p_0 (p_0')^{-1} \in R(T)$.
But, since $C(T) \cap R(T) = \{(1) \}$, we have $q_0 =q_0'$ and $p_0=p_0'$. 
Hence 
$$
v_{T}^{(0)}(1,\dots,n) =\sum_{q_0 \in C(T)} \sgn(q_0) \sum_{p_0 \in R(T)}
x_{q_0 p_0(1),1} \cdots x_{q_0 p_0(n),n}\not= 0
$$ 
and so we have proved the corollary.
\end{proof}

\begin{example}
For $T= \miniyoung(13,2)$,
we have
\begin{align*}
v_{T}^{(\alpha)}(1,2,1)
=& D^{(\alpha)} (1,2,1) \cdot ( (1)+(13) -(12)-(132)) \\
=& 2 D^{(\alpha)} (1,2,1) - D^{(\alpha)}(2,1,1) -D^{(\alpha)}(1,1,2) \\
=& (1+\alpha)(1-\alpha) (2x_{11}x_{22} x_{13} - x_{21} x_{12} x_{13} -x_{11} x_{12} x_{23}). \qed
\end{align*}
\end{example}

\begin{example}
For $T=\miniyoung(12,34)$, we have
\begin{align*}
v_{T}^{(\alpha)}(1,2,2,4) 
=& D^{(\alpha)}(1,2,2,4) + D^{(\alpha)}(1,2,4,2) -2 D^{(\alpha)}(1,4,2,2) +D^{(\alpha)}(2,1,2,4) \\
 &+ D^{(\alpha)}(2,1,4,2)  -2 D^{(\alpha)}(2,2,1,4) -2D^{(\alpha)}(2,2,4,1) + D^{(\alpha)}(2,4,1,2) \\
 &+ D^{(\alpha)}(2,4,2,1) -2D^{(\alpha)}(4,1,2,2) +D^{(\alpha)}(4,2,1,2) +D^{(\alpha)}(4,2,2,1) \\
=& (1+\alpha)(1-\alpha) ( x_{11}x_{22}x_{23}x_{44} + x_{11}x_{22}x_{43}x_{24}
  -2x_{11}x_{42}x_{23}x_{24} +  x_{21}x_{12}x_{23}x_{44} \\
  &  +x_{21}x_{12}x_{43}x_{24} -2  x_{21}x_{22}x_{13}x_{44} -2  x_{21}x_{22}x_{43}x_{14}
   +  x_{21}x_{42}x_{13}x_{24} \\
 & + x_{21}x_{42}x_{23}x_{14} -2  x_{41}x_{12}x_{23}x_{24} + x_{41}x_{22}x_{13}x_{24}
  +  x_{41}x_{22}x_{23}x_{14}). \qed
\end{align*}
\end{example}

For a semi-standard tableau $S$ and a standard tableau $T$ of the same shape,
we define the sequence $\boldsymbol{i}^{(S,T)}=(i_1^{(S,T)}, \dots, i_n^{(S,T)})$ as follows.
For each $k$, we let $B_k$ be the box numbered by $k$ in $T$ and
denote by $i_k^{(S,T)}$ the number in box $B_k$ of $S$.
For example, for
$$
S= \young(1223,334,46) \qquad \text{and} \qquad T=\young(1356,249,78),
$$
we have $\boldsymbol{i}^{(S,T)} = (1,3,2,3,2,3,4,6,4)$. 
Put 
$v_{S,T}^{(\alpha)} = v_{T}^{(\alpha)}(\boldsymbol{i}^{(S,T)})$.

\begin{example}
\begin{align*}
v_{S,T}^{(\alpha)} =&
2D^{(\alpha)}(1,2,1) - D^{(\alpha)}(2,1,1) - D^{(\alpha)}(1,1,2)  
\quad \text{for $(S,T)=\(\miniyoung(11,2),\ \miniyoung(13,2)\)$}, \\
v_{S,T}^{(\alpha)} =&
D^{(\alpha)}(1,3,2) - D^{(\alpha)}(3,1,2) + D^{(\alpha)}(2,3,1) - D^{(\alpha)}(2,1,3)
\quad \text{for $(S,T)=\(\miniyoung(12,3),\ \miniyoung(13,2)\)$}, \\
v_{S,T}^{(\alpha)} =&
D^{(\alpha)}(1,2,3) - D^{(\alpha)}(2,1,3) + D^{(\alpha)}(3,2,1) - D^{(\alpha)}(3,1,2)
\quad \text{for $(S,T)=\(\miniyoung(13,2),\ \miniyoung(13,2)\)$}.\qed
\end{align*}
\end{example}

\begin{example}
For
$$
S=T= \young(1,2,{:},n)
$$
we have $v_{S,T}^{(\alpha)}=\sum_{q \in \mf{S}_n} \sgn(q) D^{(\alpha)} (q(1), \dots, q(n))=
 \prod_{j=1}^{n-1} (1-j \alpha) \det(X)$.\qed
\end{example}

Finally, we obtain the following theorem.

\begin{thm} \label{ThmDecomposition}
Denote by $W_T^{(\alpha)}$ the image of the map $V_n^{(\alpha)} \to V_n^{(\alpha)}$
given by $c_T$
for each standard tableau $T$ of shape $\lambda \vdash n$.
Then $V_n^{(\alpha)} = \bigoplus_{\lambda \vdash n} \bigoplus_{T} W_T^{(\alpha)}$,
where $T$ run over standard tableaux,
and 
$$
W_T^{(\alpha)} \cong f_{\lambda}(\alpha) E^{\lambda} =
\begin{cases} \{ 0\} & \text{for $\alpha \in \{ 1, \frac{1}{2}, \dots, \frac{1}{\lambda_1'-1},
-1, -\frac{1}{2}, \dots, -\frac{1}{\lambda_1-1}\}$}, \\
E^{\lambda} & \text{otherwise}. 
\end{cases}
$$
When $W_T^{(\alpha)} \cong E^{\lambda}$, 
the $v_{S,T}^{(\alpha)}= f_{\lambda}(\alpha) v_{S,T}^{(0)}$,
where $S$ run over all semi-standard tableaux of shape $\lambda$ with entries in $[n]$,
form a basis  of $W_T^{(\alpha)}$.
Further, the vector
$v_{S,T}^{(\alpha)}$ is the highest weight vector of $W_T^{(\alpha)}$
if all entries in the $r$-th row of $S$ are $r$,
and $v_{S,T}^{(\alpha)}$ is the lowest weight vector if 
entries in each $r$-th column of $S$ are given as 
$n-\lambda'_r+1, \dots, n-1, n$ from the top.
\end{thm}

\begin{proof}
By Proposition \ref{PropProjection}, it is clear that 
$V_n^{(\alpha)} = \bigoplus_{\lambda} \bigoplus_{T} W_T^{(\alpha)}$
and each $W_T^{(\alpha)}$ is $\{ 0\}$ or isomorphic to $E^{\lambda}$.
The space $W_T^{(\alpha)}$ is generated by $v_{T}^{(\alpha)}(i_1,\dots, i_n)$,
where $(i_1,\dots, i_n) \in [n]^n$.
Since $W_T^{(0)} \cong E^{\lambda}$ by \eqref{EqDecompositionZero},
it follows from Corollary \ref{ThmFactorization} that
$W_T^{(\alpha)} \cong E^{\lambda}$ unless $f_\alpha(\lambda) \not=0$.
It is easy to see that
$f_{\alpha}(\lambda) =0$ if and only if 
$\alpha= 1/k$ for $ 1 \le k \le \lambda_1'-1$ or 
$\alpha=-1/k$ for $ 1 \le k \le  \lambda_1-1$.

Suppose $W_T^{(\alpha)} \cong E^{\lambda}$.
Elements $\{v_{S,T}^{(0)} \ | \ \text{$S$ are semi-standard tableaux} \}$
are linearly independent.
In fact,
for any semi-standard tableau $S_0$,
the term $D^{(0)} (\boldsymbol{i}^{(S_0,T)})= x_{i_1^{(S_0,T)},1} \cdots x_{i_n^{(S_0,T)},n}$
appears only in $v_{S_0,T}^{(0)}$ among all $v_{S,T}^{(0)}$.   
Since the dimension of $E^{\lambda}$ is equal to the number of semi-standard tableaux of shape $\lambda$,
the $v_{S,T}^{(\alpha)}=f_{\lambda}(\alpha) v_{S,T}^{(0)}$ form a basis of $W_T^{(\alpha)}$.
It is immediate to check the last claim.
\end{proof}

Theorem \ref{ThmDecomposition} says that $\{D^{(\alpha)}(i_1, \dots, i_n) \ | \ i_1, \dots, i_n \in [n]\}$
are linearly independent if $\alpha \in \bC \setminus \{\pm 1/k \ | \ k=1, \dots, n-1 \}$.
Theorem \ref{MainTheorem} follows from Theorem \ref{ThmDecomposition} immediately.

\medskip

A trick of doubling the variables (\cite{Howe}, \cite{Weyl}) suggests
the following corollary.
In fact, since $\dim_{\bC} (E^{\lambda} \otimes (E^{\lambda})^{*})^{\mf{sl}_n(\bC)} = 1$,
we can express $\det(X)^2$ by $\alpha$-determinants except a finite number of $\alpha$.

\begin{cor} \label{CorInvariants}
Let $\alpha \in \bC \setminus \{\frac{1}{k} \ | \ 1 \le k \le \frac{n-1}{2} \}$.
Then there exists $\lambda \vdash n$ such that $f_{\lambda}(\alpha) \not=0$, 
which has the following property;
for any standard tableau $T$ of shape $\lambda$
there exists a $\sl_n(\bC)$-intertwining operator $A^{(\alpha)}: (W_T^{(\alpha)})^* \to W_T^{(\alpha)}$
satisfying $A^{(\alpha)} ((v_{T,T}^{(\alpha)})^*) = v_{T,T}^{(\alpha)}$ and
\begin{equation} \label{EqInvariants}
\det(X)^2 = f_\lambda(\alpha)^{-2} 
\sum_{S} v_{S,T}^{(\alpha)} \cdot A^{(\alpha)} ((v_{S,T}^{(\alpha)})^*).
\end{equation}
Here the sum runs over all semi-standard tableaux $S$ of shape $\lambda$
and $(v_{S,T}^{(\alpha)})^*$ are defined by 
$(v_{S,T}^{(\alpha)})^* (v_{S',T}^{(\alpha)} ) =\delta_{S,S'}$.
More precisely, one may take 
$\lambda =(2^{\frac{n}{2}})$ if $n$ is even,
or $\lambda=(1^1 2^{\frac{n-1}{2}})$ if $n$ is odd,
which satisfies the condition. 
\end{cor} 

\begin{proof}
For $\alpha$ and $\lambda$ in the corollary, 
it is easy to see that $f_\lambda (\alpha) \not=0$.
Consider a standard tableau $T$ of shape $\lambda$.
Then we see that
$W_T^{(\alpha)} =W_T^{(0)}$ and $(v_{S,T}^{(\alpha)})^*=f_\lambda (\alpha)^{-1} (v_{S,T}^{(0)})^*$.
Suppose the corollary is true for $\alpha=0$.
Using the intertwining operator $A^{(0)}$,
we define $A^{(\alpha)}$ by $A^{(\alpha)} =f_\lambda(\alpha)^2 A^{(0)}$.
Then we see that 
$A^{(\alpha)} ((v_{T,T}^{(\alpha)})^*) 
=f_\lambda(\alpha)^2 A^{(0)}(f_\lambda (\alpha)^{-1} (v_{T,T}^{(0)})^*) =v_{T,T}^{(\alpha)}$ and
$$
\sum_{S} v_{S,T}^{(\alpha)} \cdot A^{(\alpha)} ((v_{S,T}^{(\alpha)})^*)
= f_\lambda (\alpha)^2 \sum_{S} v_{S,T}^{(0)} \cdot A^{(0)} ((v_{S,T}^{(0)})^*)
= f_\lambda (\alpha)^2 \det(X)^2.
$$
It is hence sufficient to prove the corollary for the case $\alpha=0$.

In general, for a finite-dimensional irreducible $U(\sl_n)$-module $V$ 
and its dual module $V^*$,
$$
\mathfrak{I}= \sum_{i} v_i \otimes v_i^* \in V \otimes V^*
$$
defines an invariant of $\sl_n(\bC)$, see \cite{Howe}.
Here $v_i$ are a basis of $V$ and $v_i^*$ are the dual basis, i.e, $v_i^*(v_j) = \delta_{ij}$. 
When $V= W_T^{(0)} \cong E^{\lambda}$ ($\lambda=(1^{n-2r} \, 2^r)$), 
the module $V$ is self-dual, i.e.,
$V^* \cong_{\sl_n(\bC)} V$. 
Therefore there exists an intertwining operator $A'$ from $V^*$ to $V$.
Then the polynomial
$\sum_{S} v_{S,T}^{(0)} \cdot A' ((v_{S,T}^{(0)})^*) 
\in \mathcal{P}(\mathrm{Mat}_{n \times n})$ of degree $2n$ determined by $\mf{I}$
is an invariant of $\sl_n$
and hence 
\begin{equation} \label{EqInvariantsConstant}
\sum_{S} v_{S,T}^{(0)} \cdot A' ((v_{S,T}^{(0)})^*) = c \det(X)^2
\end{equation}
for some constant $c$.
Comparing the coefficients of $x_{11}^2 \cdots x_{nn}^2$ in
both sides in \eqref{EqInvariantsConstant},
we have $A'((v_{T,T}^{(0)})^*) = c x_{11} \cdots x_{nn} =c v_{T,T}^{(0)}$
and so $c \not=0$.
Hence $A^{(0)}= c^{-1} A'$ is our desired operator.
\end{proof}

\begin{example}
Let $T=\miniyoung(12)$. The module $W_T^{(\alpha)}$ has 
a basis consisting of 
$v_+=v_{\mminiyoung(11),\mminiyoung(12)}^{(\alpha)} =2 D^{(\alpha)}(1,1)$, 
$v=v_{\mminiyoung(12),\mminiyoung(12)}^{(\alpha)} =D^{(\alpha)}(1,2) + D^{(\alpha)}(2,1)$,
$v_-=v_{\mminiyoung(22),\mminiyoung(12)}^{(\alpha)} =2 D^{(\alpha)}(2,2)$
if $\alpha \not=-1$.
The linear map $A$ determined by 
$$
A (v_+^*)
= -\frac{1}{2} v_-, \quad 
A (v^*)
= v, \quad
A (v_-^*)
= -\frac{1}{2} v_+
$$
from $(W_T^{(\alpha)})^*$ to $W_T^{(\alpha)}$
defines an intertwining operator of $\sl_2(\bC)$.
Hence, by the corollary,
we have
\begin{align*}
(1+\alpha)^2 \det(X)^2 =&
v_+ \cdot A (v_+^*)
+ v \cdot A (v^*)
+v_- \cdot A (v_-^*) 
= v^2 - 
v_+ \cdot v_- \\
=& (D^{(\alpha)}(1,2) + D^{(\alpha)}(2,1))^2 -4 D^{(\alpha)}(1,1)D^{(\alpha)}(2,2). \qed
\end{align*}
\end{example}

%
\section{Concluding remarks}
%

We give here a sketch of another proof of the main theorem based
on a character theory of $\mf{S}_n$ together with a comment 
on a possible generalization of the theorem 
to the quantum group $U_q(\gl_n)$ 
(Section \ref{SectionQuantum} and \ref{SectionImmanant}).
At the final position of the paper, 
we discuss the irreducible decomposition of $V_n^{(\alpha)}$
when $\alpha=\infty$ (Section \ref{SectionInfty}).

\subsection{}\label{SectionQuantum}

Using representation theory of the symmetric groups,
a simpler proof of Proposition \ref{ThmCycleFormula} can be given as follows. 
The formula
\begin{equation} \label{EqFCF}
\alpha^{n-\nu_n(\sigma)} = \sum_{\lambda \vdash n} 
\frac{f^\lambda}{n!} f_\lambda(\alpha) \chi^\lambda(\sigma)
\end{equation}
is a specialization of the Frobenius character formula, see \cite[I-7, Example 17]{Mac}.
Here $\chi^\lambda$ is the irreducible character of $\mf{S}_n$ corresponding to $\lambda$.
Moreover, 
for a standard tableau $T$ and a partition $\lambda$, 
the well-known equation 
\begin{equation} \label{EqYoung}
\chi^\lambda \cdot c_T = \delta_{\lambda, \mathrm{sh}(T)} \frac{n!}{f^\lambda} c_T
\end{equation}
in the group algebra $\bC\mf{S}_n$ 
holds, which is obtained by Young.
Here $\mathrm{sh}(T)$ is the shape of $T$. 
Proposition \ref{ThmCycleFormula} follows from \eqref{EqFCF} and \eqref{EqYoung}.

However, 
when we develop our discussion to the quantum enveloping algebra,
we know that the $q$-Young symmetrizer is not essential.
Define the quantum $\alpha$-determinant by
\begin{equation} \label{DefQalpha}
\dete{\alpha,q}(X) = \sum_{\sigma \in \mf{S}_n} q^{\mathrm{inv}(\sigma)}
\alpha^{n-\nu(\sigma)} x_{\sigma(1)1} \cdots x_{\sigma(n) n},
\end{equation}
where $\mathrm{inv}(\sigma)$ is the inversion number of $\sigma$;
$\mathrm{inv}(\sigma) =\# \{ (i,j) \ | \ 1 \le i< j \le n, \ \sigma(i) > \sigma(j) \}$.
In particular, $\dete{q}(X)=\dete{-1,q}(X)$ is the 
(usual) quantum determinant, 
see e.g. \cite{NUW}.
We define quantum analogue elements 
$v_{S,T}^{(\alpha,q)}$  of $v_{S,T}^{(\alpha)}$
by the $q$-Young symmetrizers studied by Gyoja \cite{Gyoja}.
For each standard tableau $T$, denote by $W_T^{(\alpha,q)}$
the quantum analogue of $W_T^{(\alpha)}$
given by the $q$-Young symmetrizer.
Let $\lambda$ be a partition of $n$ and let
$T_1, \dots, T_d$ be all standard tableaux of shape $\lambda$,
where $d=f^\lambda$.
Let $v_{k}^{(\alpha,q)}=v_{S_k,T_k}^{(\alpha,q)}$ be the highest weight vector 
of each $W_{T_k}^{(\alpha,q)}$.
Then
there exists a $d \times d$ matrix $F_\lambda (\alpha;q)$ such that
$$
(v_{1}^{(\alpha,q)}, \dots, v_{d}^{(\alpha,q)})
= (v_{1}^{(0,q)}, \dots, v_{d}^{(0,q)}) F_\lambda(\alpha;q).
$$
In the classical case, as we have seen in Corollary \ref{ThmFactorization},
$F_\lambda(\alpha;1)$ is the scalar matrix $f_{\lambda}(\alpha) I$.
It is observed, however, 
$F_\lambda(\alpha;q)$ is not, in general, 
a scalar matrix, not even a diagonal matrix,
see \cite{KMW}.
Therefore, it is necessary to find a new bases other than 
the one obtained by the $q$-Young symmetrizers in order to 
diagonalize $F_\lambda(\alpha;q)$.

\subsection{}\label{SectionImmanant}

Recall the immanant.
For a partition $\lambda$ of $n$, the $\lambda$-immanant of $X$ is defined by
$$
\imm(X)= \sum_{\sigma \in \mf{S}_n} \chi^\lambda (\sigma) \prod_{i=1}^n x_{i \sigma(i)}.
$$
Then, if we use the formula \eqref{EqYoung}, we find 
the cyclic module $U(\gl_n) \imm(X)$ is decomposed as 
$U(\gl_n) \imm(X) \cong (E^{\lambda})^{\oplus f^\lambda}$
as in the case of $\alpha$-determinants. Also, 
since the function $\sigma \to \nu_n(\sigma)$ is a class function,
the $\alpha$-determinant is expanded by immanants;
$$
\dete{\alpha}(X) = \sum_{\lambda \vdash n} 
\frac{f^\lambda}{n!} f_{\lambda}(\alpha) \imm(X)
$$
by \eqref{EqFCF}.
Combining these facts, we have
$V_n^{(\alpha)} \cong 
\bigoplus_{\begin{subarray}{c} \lambda \vdash n \\ f_{\lambda}(\alpha) \not= 0 \end{subarray}}
(E^{\lambda})^{\oplus f^\lambda}$.
This agrees with Theorem \ref{MainTheorem}.
However, if we consider the quantum case,
this discussion can not be applied,
because the function $\sigma \mapsto \mathrm{inv}(\sigma)$ 
appeared in \eqref{DefQalpha}
is not a class function.

\subsection{}\label{SectionInfty}

We consider the case ``$\alpha=\infty$''.
Since $\dete{\alpha}(X)$ is 
a polynomial of degree $n-1$ in variable $\alpha$,
we can define a limit
\begin{equation} \label{DefDetInfty}
\dete{\infty} (X)= \lim_{|\alpha| \to \infty} \alpha^{1-n} \dete{\alpha}(X)
= \sum_{\begin{subarray}{c} \sigma \in \mf{S}_n \\ \nu_n(\sigma)=1 \end{subarray}} 
x_{\sigma(1),1} \cdots x_{\sigma(n),n}.
\end{equation}
For example, 
$$
\dete{\infty} \begin{pmatrix} x_{11} & x_{12} & x_{13} \\ x_{21} & x_{22} & x_{23}
\\ x_{31} & x_{32} & x_{33} \end{pmatrix}
= x_{21} x_{32} x_{13} + x_{31} x_{12} x_{23}.
$$
Denote by $V_n^{(\infty)}$ the cyclic module $U(\mf{g}) \dete{\infty}(X)$.
Then we have the following irreducible decomposition of $V_n^{(\infty)}$ 
as a corollary of Theorem \ref{ThmDecomposition}.

\begin{cor}
$$
V_n^{(\infty)} \cong \bigoplus_{\lambda : \mathrm{hook}} (E^\lambda)^{\oplus f^\lambda}
= \bigoplus_{k=1}^n \(E^{(k,1^{n-k})} \)^{\oplus \binom{n-1}{k-1}},
$$
where $\lambda$ run over all hook partitions of $n$.
\end{cor}

\begin{proof}
The degree of polynomial $f_{\lambda}(\alpha) \in \bZ[\alpha]$ is equal to $n -d$, 
where $d$ is the number of the main diagonal of the Young diagram $\lambda$.
Therefore $\lim_{|\alpha| \to \infty} \alpha^{1-n} f_{\lambda}(\alpha)$ is zero
unless $d=1$, i.e., $\lambda$ is a hook.
For a hook $\lambda=(k,1^{n-k})$, the number $f^\lambda$ is given 
by the binomial coefficient $\binom{n-1}{k-1}$.
Hence, the claim follows from Theorem \ref{ThmDecomposition}.
\end{proof}

\begin{example}
When $n=5$,
$$
V_5^{(\infty)} \cong E^{(5)} \oplus \( E^{(4,1)} \)^{\oplus 4} \oplus \( E^{(3,1,1)} \)^{\oplus 6}
\oplus \(E^{(2,1,1,1)} \)^{\oplus 4} \oplus E^{(1,1,1,1,1)}.  \qed
$$
\end{example}

\begin{remark}
By \eqref{EqStanley}, for each $\alpha > 0$,  
we can define a probability measure $\mathfrak{M}_n^{(\alpha)}$ on $\mf{S}_n$ by
$$
\mathfrak{M}_n^{(\alpha)} (\sigma)= \frac{\alpha^{n-\nu_n(\sigma)}}{\prod_{j=1}^{n-1} (1+j \alpha)}
\qquad \text{for each $\sigma \in \mf{S}_n$}.
$$
This is called the Ewens measure in \cite{KOV} 
but the definition is slightly different from ours. 
It is clear that
$\mathfrak{M}_n^{(1)}$ is the uniform measure on $\mf{S}_n$
and $\mathfrak{M}_n^{(0)} = \lim_{\alpha \to 0+} \mathfrak{M}_n^{(\alpha)}$
is the Dirac measure at the identity.
Also we see that
$$
\mathfrak{M}_n^{(\infty)} (\sigma)
= \lim_{\alpha \to +\infty} \mathfrak{M}_n^{(\alpha)} (\sigma)
= \begin{cases} 1/(n-1)! & \text{if $\nu_n(\sigma)=1$}, \\
 0& \text{otherwise}. \end{cases}
$$
Given a matrix $X$ with non-negative entries $x_{ij}$,
define a random variable $X_\sigma$ by $X_{\sigma} =\prod_{i=1}^n x_{i\sigma(i)}$ on $\mf{S}_n$.
Then for $\alpha \in [0, +\infty]$ the $\alpha$-determinant of $X$ is essentially 
the mean value of $X_\sigma$ with respect to $\mathfrak{M}_n^{(\alpha)}$:
\begin{align*}
\dete{\alpha}(X) 
=& \prod_{j=1}^{n-1} (1+j \alpha) 
\sum_{\sigma \in \mf{S}_n} X_\sigma \mathfrak{M}_n^{(\alpha)}(\sigma) 
\qquad \text{for $0 \le \alpha < +\infty$}, \\
\dete{\infty}(X) 
=& (n-1)! \sum_{\sigma \in \mf{S}_n} X_{\sigma} \mathfrak{M}_n^{(\infty)}(\sigma). \qed
\end{align*}
\end{remark}

\vspace{11pt}

\noindent
{\bf Acknowledgement.} 
We would like to thank Kazufumi Kimoto for his helpful comments.


\noindent
\textsc{Sho MATSUMOTO}\\
Graduate School of Mathematics, Kyushu University.\\
Hakozaki Higashi-ku, Fukuoka, 812-8581 JAPAN.\\
\texttt{shom@math.kyushu-u.ac.jp}\\

\noindent
\textsc{Masato WAKAYAMA}\\
Faculty of Mathematics, Kyushu University.\\
Hakozaki Higashi-ku, Fukuoka, 812-8581 JAPAN.\\
\texttt{wakayama@math.kyushu-u.ac.jp}

\end{document}